\documentstyle[a4]{article}
\newtheorem{thm}{Theorem}
\newtheorem{prop}{Proposition}
\newtheorem{cor}{Corollary}
\newtheorem{lem}{Lemma}
\newcommand{\tree}[1]{\mbox{\scriptsize $\left(#1 \right)$}}
\newcommand{\nin}{\not\in}

\begin{document}
\begin{titlepage}
\title{A New Proof of Monjardet's Median Theorem}
\author{D. E. Loeb\thanks{Author supported by URA CNRS 1304 and
NATO CRG 930554.}\\ 
LaBRI\\
Universit\'e de Bordeaux I\\
33405 Talence, France\\
loeb@labri.u-bordeaux.fr}
\end{titlepage}
\maketitle

\begin{center}
Dedicated to Bernard Monjardet
\end{center}

\begin{abstract}
New proofs are given for Monjardet's theorem that all strong simple
games (i.e., ipsodual elements of the free distributive lattice) can
be generated by the median operation. Tighter limits are placed on the
number of iterations necessary. Comparison is drawn with the $\chi$
function which also generates all strong simple games.
\end{abstract}

\section{Introduction}
In the preceding article, \cite{L} we studied the composition of
several games. In such a composition, the voters of one game are
replaced by committees each voting according to the rules of the other
games. Here, we remove the condition that the committees be disjoint. 
In fact, adding dummy voters when necessary, we may assume without
loss of generality, that each of the committees is in fact the entire
set of voters. Thus, in the generalized compound games considered
here, the $n$ voters vote $k$ times according to the rules
of various games. These $k$ results are combined via a $k$-voter
game. This process defines a new
game, preserving the {\em strength} and/or 
{\em simplicity} of the original games.

As opposed to the situation with disjoint committees, here there is no
sense any {\em irreducible} strong simple games. Nevertheless, the 
three-voter 
democracy $Dem_3$ together with the $n$ dictatorships can be singled out as
{\em basic} games, since when combined using this generalized
composition, they generate all $n$-voter strong simple games.
This fact was proven by Monjardet \cite{M1,M2}. However, a new proof
is given here which makes no reference to games which are not both simple
and strong. This new proof sets a tighter bound on the number of
compositions required to generate all strong simple games, yet is
completely elementary requiring no lattice theory. 

To generate all $n$-voter strong games, an additional
basic game 
is needed. It is the  winning  game
$$\hat{1} =\{ \mbox{all subsets of $V$}
\} = (0\dots 0)_0.$$
The notation
$[a_{1}a_{2}\cdots ]_q$ means that each voter $i$ is given $a_{i}$
``votes.'' A coalition wins if it attains the quota of
$q$ votes. Any 
game which can be expressed in this notation is called a {\em quota}
game.

Similarly, to generate all $n$-voter simple games, the 
losing game $\hat{0}$ is needed. 
$$\hat{0}=\emptyset = (0\dots 0)_1.$$ 
To generate all games, both $\hat{0}$ and $\hat{1}$ are needed as
additional basic games.

Monjardet's proof makes use of the well-known identification of 
$n$-voter games with points in the free-distributive lattice on
$n$-generators $FD_{n}$.  \cite{S}
The generators $i$ themselves represent {\em dictatorial games} in
which $C$ is a winning coalition if and only if the dictator $i\in
C$. 

Given three elements of a distributive lattice $S,T,U$, one defines
their {\em median} to be the meet of their joins:
$$ m(S,T,U)=ST+SU+TU=(S+T)(S+U)(T+U). $$
In terms of games, the median law is the composition of the three
games with $Dem_3$.

$FD_{n}$ is a self-dual lattice. (The dual of any given element
can be found by expressing it in terms of the generators and
exchanging ``meet'' with ``join'' throughout.) Many interesting
properties of games can be expressed in terms of this duality. 
In a strong game $S$ (resp. simple, resp. strong simple), we have
$S\geq S^*$ in $FD_n$ (resp. $S\leq S^*$, resp. $S=S^*$). In the
language of free
distributive lattices, this property is called {\em supraduality}
(resp. {\em infraduality}, resp. {\em ipsoduality}).complement loses. 

It is obvious that if $S,T,U\in FD_{n}$ are all ipsodual (resp.
infradual, resp. 
supradual), then so will be their median, since
$$ m(S,T,U)^{*} = m(S^{*},T^{*},U^{*}).$$ 
However, less obvious is the
result proven by Monjardet that all ipsodual elements can
be generated from the $n$ dictatorial games via repeated application
of the median operation. 
\begin{thm}[Monjardet \cite{M1,M2}] \label{monj}
Let $FM_{n}$ be the smallest set containing the generators of $FD_{n}$
and closed under the median operation. Given $S\in FD_{n}$, $S$ is
ipsodual if and only if $S\in FM_{n}$.
\end{thm}
In other words, strong simple games can be identified with elements of
the free median set. That is to say, all games can be
thought of as compositions of  triples of simpler games,
where the simplest games of all are dictatorial.

See \cite{L} for additional motivation and notation.

\section{Quotient Games}
Let $S$ be a strong simple game with player set $V$, and let
$f:V\rightarrow W$ 
be some function. Define the {\em quotient game}
$f(S)=\{A\subseteq W: f^{(-1)}(A)\in S \}$. 
$V$ can be thought of as a set of {\em offices} and $W$ as a set of
{\em voters}. $f$ describes which offices are held by which voters. 

Voters which appear in no minimal winning coalitions are called {\em
dummies}, otherwise voters are called {\em powerful}.
If
$f$ is non-surjective, then certain voters will hold no office, and
thus be {\em dummies}. If $f$ is non-injective, then certain
voters will 
combine the functions of several offices. A single vote of a voter is taken
into account as the vote of each of his offices. If $f$ is bijective, then
$f(S)$ is isomorphic to $S$.

Note that the minimal
winning coalitions of $f(S)$ are the images of certain minimal
winning coalitions of $S.$ 

If $S$ is a quota game, then $f(S)$ will also be a quota game in which
the weight of each voter is equal to the sum of the weights of his
offices.  

Our first proof of Theorem 1 depends on the following two lemmata.
\begin{lem}
Let $S$ be a game with voter set of $V$, and let
$f:V\rightarrow W$ 
be some function. Then $f(S)$ is a game on $W.$ If $S$ is strong
(resp. simple), then $f(S)$ is too.$\Box $
\end{lem}

{\em Proof:} $f^{(-1)}$ is a monotone function from $2^{W}$ to $2^{V}$,
thus $f(S)$ is a game. Suppose $S$ is simple, and $B\in f(S)$. Then
$f^{(-1)}(B)\in S$. Hence, $V-f^{(-1)}(B)=f^{(-1)}(W-B)\nin S $.
Strength follows as $W-B\nin f(S).\Box $

\begin{lem}\label{fig}
All nondictatorial strong simple games are the medians of strong
simple games with a strictly greater number of dummies.
\end{lem}

{\em Proof:} Let $x,y,z$ be distinct powerful voters. Let
$f_{ij}:V\rightarrow 
V$ fix all members of $V$ with the exception of $i$ for which
$f_{ij}(i)=j$. We will show that 
\begin{equation}\label{med}
S=m(f_{xy}(S),f_{yz}(S),f_{zx}(S)). 
\end{equation}
If so, then we have expressed $S$ in terms of strong simple games
involving at most $n-1$ powerful voters. 

Let $T$ denote the right hand side of equation (\ref{med}). 
Since $|S|=|T|=2^{n-1}$, it will suffice to show that $S\subseteq T$.
Suppose, $A\in S$. Consider then the possible values of $|A\cap \{x,y,z
\}|:$
\begin{itemize}
\item [0.] If $x,y,z\nin A$. Then $A\in f_{xy}(S),f_{yz}(S),f_{zx}(S)$.
Thus, $A\in T$.
\item [1.] Without loss of generality, suppose $x\in A$, and $y,z\nin
A$. Then $A\in f_{yz}(S),f_{zx}(S)$, so $A\in T$.
\item [2.] Without loss of generality, suppose $x,y\in A$, and $z\nin
A$. Then $A\in f_{xy}(S),f_{zx}(S)$, so $A\in T$.
\item [3.] If $x,y,z\in A$. Then  $A\in f_{xy}(S),f_{yz}(S),f_{zx}(S)$.
Thus, $A\in T.\Box $
\end{itemize}

{\em Proof of Theorem 1:} Since the number of dummies is bounded by
$n-1$, Monjardet's theorem now follows by induction. QED.

We have the following result more generally.

\begin{thm}
Let $S$ be a game (resp. strong game, resp. simple game, resp. strong
simple game) with at least two powerful voters. Then $S$ is
the median of three games (resp. strong games, resp. simple games, resp. strong
simple games) with a strictly greater number of dummies.
\end{thm}

{\em Proof:} If $S$ is a strong simple game, then the lemma above
suffices since all nondictatorial games have at least two powerful
voters.

If $S$ has three or more powerful voters, then the above reasoning is still
valid. 

If $S$ has exactly two powerful voters (say $x$ and $y$) then $S$ must
either be the simple majority game $xy=(1,1,0,0,\dots)_2$ or the strong majority
game $x+y=(1,1,0,0,\dots)_1$. In the game $xy$, the participation of both
$x$ and $y$ are needed to win, whereas in $x+y$ the participation of
either is sufficient. Notice that $xy=m(Dict_x,Dict_y,\hat{1})$ and
$x+y=m(Dict_x,Dict_y,\hat{0})$ where $Dict_x$ and $Dict_y$ denote
dictatorship by $x$ and $y$ respectively. (Note that in the games
$\hat{0}$ and $\hat{1}$, all voters are dummies.$\Box$

\begin{cor}
All games (resp. strong games, resp. simple games) can be generated
via generalized composition of various 
dictatorships, and the three-voter democracy $Dem_3=(1,1,1)_2$
along with
$\hat{0}$ and $\hat{1}$
(resp. $\hat{0}$ alone, resp. $\hat{1}$ alone).
\end{cor}

{\em Proof:} By induction, all games can be reduced via $Dem_3$ to games of the
same ``type'' having only 0 or 1 powerful voters. The only games with
1 powerful voter are the dictatorships. The only games with no
powerful voters are the predetermined win $\hat{1}$ and the
predetermined loss $\hat{0}$. The former is strong and the latter is
simple.$\Box$
 
Note that the dictatorships along with $\hat{0}$, $\hat{1}$,
$or=(1,1)_1$, and $and=(1,1)_2$ also forms a minimal basis from which
all games can be constructed. In other words, all monotonic boolean
functions can be written as a for example conjunction of disjunctions. 
 
\section{Weight}
Define the weight of a strong simple game $S$, to be
the number of times the 
median operation must be iterated in order
to generate $S$ from dictatorships. 
In other words, a
game $S$ is of weight
$0$ if it is a dictatorship. For $k>0$, $S$ is of 
weight $k$ if it is the median of three strong simple
games of weight $k-1$ or less, but not the median of three strong
simple games of weight $k-2$ or less.

Define $Dem_3^0$ to be the unique one-voter strong simple game $Dem_1$, and
$Dem_3^{d+1}=Dem_3[Dem_3^d,Dem_3^d,Dem_3^d]$. 
$Dem_3^d$ is clearly a transitive strong simple game on $3^d$ voters. 
A strong simple game has weight at most $d$ if and only if it is
quotient of $Dem_3^d$.

Let $W(n)$ be the largest weight of an
$n$-voter strong simple game. 

\begin{prop} (1) For $n\leq 6$, $W(n)$ is given by the following table.
$$ \begin{array}{r|rrrrrrrr}
n &   1 & 2 & 3 &  4 &  5 & 6  \\ \hline
W(n)& 0 & 0 & 1 &  2 &  3 & 3 
\end{array} $$

(2) For $n\geq 6$, $W(n)\leq n-3$.

(3) $W(n)$ is asymptotically greater than $\frac{\ln 2}{\ln 3}n$.

(4) $W(n)$ is weakly increasing, restricted growth function. ie. $W(n+1)$
always equals either $W(n)$ or $W(n)+1$.
\end{prop}

{\em Proof:} (1) See above.

(2) By the proof of theorem \ref{monj} above, the weight of a strong
simple game on $n$ voters is no more than one greater than the largest
such weight on $n-1$ voters.

(3) Every $n$-voter strong simple game is a quotient of
$Dem_3^{W(n)}$. There are 
$n^{3^{W(n)}}$ such quotients. Whereas, (see \cite{K}) for $n$ odd,
the number of strong simple games is asymptotically
$$ 2^{n-1 \choose (n-1)/2}\exp \left(
\parbox{3.15in}{$\displaystyle {n-1 \choose (n-1)/2}
\left(2^{-(n-1)/2} + 3n^2\dot 2^{-n-4} - n \dot 2^{-n-2} \right) + $\\
$\displaystyle {n-1 \choose (n+3)/2} \left( 2^{-(n+3)/2} + n^2 \dot
2^{-n-6} - n \dot 2^{-n-5}\right)$} \right), $$
and for $n$ even,
$$ 2^{n-1 \choose n/2}\exp \left(
\parbox{3.15in}{$\displaystyle {n-1 \choose n/2-1}
\left(2^{-n/2-1} + n \dot 2^{-n-4} \right) + $\\
$\displaystyle {n \choose n/2+1} \left( 2^{-n/2-1} + n^2 \dot
2^{-n-5} - n \dot 2^{-n-4}\right)$} \right). $$
In either case, taking the double $\ln$ of these asymptotic formulas
yields 
$$ W(n) \succeq {\frac {n\ln 2}{\ln 3}}-\frac{\ln \ln n}{\ln
3}-\frac{\ln \left( (\ln 2)^2 / 2\pi \right)}{\ln 9} + O(n^{-1}).$$

(4) $W(n)$ is a restricted growth function by the proof of theorem
\ref{monj}. Now, let $S$ be an $n$-voter
strong simple game of weight $W(n)$. Let $\iota(S)=S\cup \{A\cup
\{n+1 \}:A\in S \}$ be the corresponding $n+1$-voter game in which the
additional voter is powerless. By hypothesis, $\iota(S)$ is a quotient
of $Dem_{3}^{W(n+1)}$. However, $S=\iota(S)_{n,n+1}$ is a quotient of
$\iota(S)$. Thus, $S$ is a quotient of $Dem_{3}^{W(n+1)}$. Hence,
$W(n+1)\geq W(n).\Box $ 

\section{Choice Function}
A finite multi-player deterministic sequential-move perfect-knowledge
game can be represented by a labeled rooted tree. The 
label of an internal node indicates which player must move. His move
consists of the selection of a child of that node. The label of a leaf
indicates the winner of the game. Without loss of generality, we can consider
only choice binary trees, since choices between a large
number of possibilities can be made via iterated binary choices.

The win-type of a game-tree is the set of winning coalitions. A
coalition is said to be winning if it has a combined strategy that
ensures that the winner of the game will be one of its members.

Let $\tau$ be the game-tree 
\begin{equation}\label{choice}
\tree{\begin{array}[c]{c@{}c@{}c@{}c@{}c}
&& a\\[-.5mm]
&\swarrow&&\searrow\\[-1mm]
\tau_{1}&&&&\tau_{2}
\end{array}},
\end{equation}
where $\tau_{1}$ and $\tau_{2}$ are sub-game-trees with win-types
$S_{1}$ and $S_{2}$ respectively, and $a\in V$ is a player.
What is the win-type $S$ of the game $\tau$? 

Obviously, $A\subseteq V$ wins $\tau$ if it wins both $\tau_{1}$ and
$\tau_{2}$. 
Moreover, if $a\in A \subseteq V$, then $A$ can win $\tau$ even if it
wins only one of $\tau_{1}$ and $\tau_{2}$.
Thus, $S=({S_{1}}\cap {S_{2}} \cup \{A\in
S_{1} \cup S_{2} : a\in A \}$ 
We will denote this combination of $S_{1}$ and $S_{2}$ by 
$ \chi_{a}(S_{1},S_{2}) $ and call it the choice by $a$ between
$S_{1}$ and $S_{2}$.

\begin{prop}\label{win}
(1) The win-type of a game-tree is a strong simple game.

(2) Conversely, all strong simple games are the win-type of some
game-tree. 
\end{prop}

{\em Proof:} {\bf (1)} By induction on game-trees.

If the height is zero, then the tree is trivial. It contains one node.
The game is an automatic win for some player. A coalition is winning
if and only if it contains this player. The win-type is a dictatorship
by this player.

Otherwise, $\tau$ is of the form $\chi_{a}(S_{1},S_{2})$ where $S_{1}$
and $S_{2}$ are strong simple games. Let $a\in A$ and $B=X\setminus
A$. If $A\nin 
\chi_{a}(S_1,S_2)$, then $A\nin S_1$ and $A\nin {\cal
G}$. Thus, $B\in S_1$ and $B\in S_2$. Thus, $B\in
\chi_{a}(S_1,S_2)$. Conversely, if $A\in \chi_{a}({\cal
F},S_2)$, then $A\in S_1$ or $A\in S_2$. In either
case, $B\nin \chi_{a}(S_1,S_2).$  

{\bf (2)} Consider the following game. Players take turns
eliminating all but one of the coalitions in $S$ to which they
belong (if any). 
After each player has taken his turn there will be one set left. The
first member of this set wins.$\Box$ 

Thus, the functions $\chi_{a}$ together with the dictatorships can be
used to generate all strong simple games. 

{\em Alternate Proof of Theorem 1:} Notice that
$\chi_{a}(S,T)=m(S,T,Dict_{a})$; 
\begin{eqnarray*}
m(S,T,Dict_{a}) &=& (S\cap T)\cup (S\cap Dict_{a})\cup (T\cap Dict_{a})\\
&=& (S\cap T)\cup \{A:a\in A\in S \}\cup \{A:a\in A\in T \}\\
&=& \chi_{a}(S,T).
\end{eqnarray*}
Since $\chi$ generates all strong simple games, a fortiori $m$ does.$\Box $

\section{Depth}
In analogy to Monjardet's definition of weight, 
define the  depth of a strong simple game $S$, to be
the number of times the 
choice functions must be iterated in order
to generate $S$ from dictatorships. 
In other words, a game $S$ is of depth
$0$ if it is a dictatorship. 
$S$ is of depth $k>0$ if it is the choice by some player between two
strong simple 
games of depth $k-1$ or less, but not between two strong
simple games of depth $k-2$ or less.

$S$ is of depth $k$ if it is the win-type of a game-tree of height $k$
and no such game-tree of lesser height. Let $B_{k}$ be the complete
binary tree of height $k$. Label the nodes of $B_{k}$ consecutively
with the integers 1 to $2^{n+1}-1$. 
$$ B_{0}=(1),\mbox{ }B_{1}= \tree{\begin{array}[c]{r@{}c@{}l}
&1\\[-.5mm]
&\swarrow\searrow\\[-.5mm]
2&&3\end{array}}\mbox{, and }
B_{2}=
\tree{\begin{array}[c]{r@{}r@{}c@{}l@{}l}
&&1\\[-.5mm]
&&\swarrow\searrow\\[-.5mm]
&2&&3\\[-.5mm]
&\swarrow \downarrow && \downarrow \searrow \\[-.5mm]
4&5&&6&7
\end{array}} $$
$S$ is of depth $k$ if and only if
it is a quotient of the $2^{n+1}-1$ voter strong simple game $B_{k}$. 

Let $D(n)$ be the largest depth of an $n$-voter strong simple game. 

\begin{prop} (1) The depth of any strong simple game is at least as
great as its weight. Thus, $D(n)\geq W(n).$

(2) $D(n)$ is nondecreasing.

(3) $D(n)$ is asymptotically greater than $n$.

(4) $D(n)$  is asymptotically less than $\pi^{2}n^{3}/72$. 

(5) Let $S$ be a strong simple game of weight $w$. $S$ has depth at
most $2^w-1$. This bound can not be improved in general.
\end{prop}

{\em Proof:} {\bf (1)} $\chi_{a}(S,T)=m(S,T,Dict_{a})$. 

{\bf (2)} Let $S$ be an $n$-voter strong simple game of depth $D(n)$.
Let $\iota(S)$ be the corresponding $n+1$-voter game in which the
additional voter is powerless. Let $\iota(S)$ be the win-type
of an $n+1$-player game-tree of depth $d$. Since the additional player
is powerless, we can arbitrarily relabel nodes of the game-tree to
form an $n$-player game-tree of win-type $S$. Thus, $D(n+1)\leq d\leq D(n)$.

{\bf (3)} There are $n^{2^{d+1}-1}$ ways to label $B_{d}$ with $n$
players. Thus, there are at most $n^{2^{d+1}-1}$ $n$-player strong
simple games of depth $d$ or less. Compare with \cite{K}.

{\bf (4)} Proof of theorem \ref{win}.

{\bf (5)} Let $\tau$ (resp. $\tau',\tau'')$ be the game-tree of
$Dem_{3}^{w-1}$ with voters 
labelled $1$ to $3^{w-1}$ (resp. $3^{w-1}+1$ to $2\times 3^{w-1}$,
$2\times 3^{w-1}$ to $3^{w}$). Place a copy of $\tau'$ and $\tau''$
under each leaf of $\tau$. The resulting tree is the game-tree of
$Dem_{3}^{w}$. By induction, the height of the tree is $2^{w}-1$. 

No smaller tree can represent $Dem_{3}^{w}$ since every vertical chain
in a game-tree gives rise to a winning coalition and $Dem_{3}^{w}$ has
no winning coalitions smaller than $2^{w}.\Box $

Given a strong simple game $S$, and two voters $x,y\in V$. Say that
$x$ is at least as {\em influential} as $y$, written $x\geq y$ or
$x\geq_{S} y$ if for all
$A\in S$, we have $(A-\{y \})\cup \{x \}\in S$. 

\begin{prop}
Influence is a pre-order on the set of voters. 
\end{prop}

{\em Proof:} {\bf (Reflexivity)} Let $x\in A\in S$. Then $(A-\{x \})
\cup \{x \}=A\in S$. 

{\bf (Transitivity} Let $x\leq y\leq z$. Let $x\in A\in S$. Then $y\in
B=(A-\{x \})\cup \{y \}\in S$. Hence, $C=(B-\{y \})\cup \{z \}\in S$.
Note that $C=(A-\{x \})\cup \{z \}.\Box $

If $S$ is a quota game, then the influence relation is total.
$w_{x}\geq w_{y}$ implies $x\geq y$.\footnote{Nevertheless, the
influence relation on $S_{6,23}$ is total 
$1>2\sim 3> 4\sim 5>6$, yet $S_{6,23}$ is not a quota game. See \S{A}.}

In the proof of lemma \ref{fig}, we defined $f_{xy}(S)$. Informally,
$S_{xy}=f_{xy}(S)$ is the voting scheme in which $x$ 
``leaves the room'' having left instructions to vote according to $y$.
Thus, 
$$ S_{xy} = \{A : \mbox{$y\nin A$ and $A\in S$, or $y\in A$ and $A\cup
\{x \}\in S$} \}. $$ 
Similarly, we can imagine a voting scheme $S_{x/y}$ in which $x$
``leaves the room'' having left instructions to vote {\em against}
$y$. Thus,
$$ S_{x/y} = \{A :\mbox{$y\in A\in S$, or $y\nin A$ and $A\cup \{x
\}\in S$} \} .$$ 
\begin{prop}
Let $S$ be a game. $S_{x/y}$ is a game if and only if $x$ is less
influential than $y$, $x\leq _{S} y$.
\end{prop}

{\em Proof:} $S_{x/y}$ is a game if and only if whenever $A\in
S_{x/y}$, we have $A\cup \{a \}\in S_{x/y}$  for all $a\in V$. 

Suppose $A\in S_{x/y}$. If $a\neq y$, then $A\cup \{a \}\in
S_{x/y}$, since $S$ is a game.

On the other hand, let $a=y$. The condition $A\cup \{y \}\in
S_{x/y}$ is of interest only if $y\nin A$. In that case, it reduces to
the question of whether $A\cup \{x \}\in S$ implies $A\cup \{y \}\in
S$. In other words, whether $x\leq_{S} y.\Box $

Let $S=(w_{1},w_{2},\ldots,w_{n})_{q}$ be a quota game. Without loss of
generality, suppose $w_{1}\geq w_{n}$. Then
\begin{eqnarray}\label{part1}
S_{n1}&=&(w_{1}+w_{n},w_{2},w_{3},\ldots,w_{n-1},0)_{q}\mbox{ and}\\
\label{part2}
S_{n/1}&=&(w_{1}-w_{n},w_{2},w_{3},\ldots,w_{n-1},0)_{q-w_{n}}
\end{eqnarray}
are both quota games. If $S$ was a {\em homogeneous} or {\em
efficient} quota game (i.e., all minimal 
winning coalitions equal the total), then so is $S_{n1}$ and $S_{n/1}$.

The opposition $S_{x/y}$ plays the same role in choice function $\chi$
decompositions as the quotient $S_{xy}=f_{xy}(S)$ played in median
decompositions (lemma \ref{fig}).

\begin{prop}\label{decompose}
 Let $S$ be a strong simple game, and suppose $i\leq_{S} j$. Then 
$$ S=\chi_{i}(S_{ij},S_{i/j}). $$
\end{prop}

{\em Proof:} Say $A$ {\em needs} $i$ if $A\in S$ and $A-\{i\}\nin S.$  
Clearly, if $A\in S$ does not need $i$, then $A\in S_{ij}$ and $A\in
S_{i/j}$. Thus, $A\in \chi_{i}(S_{ij},S_{i/j})$. 

If $A\in S$ needs $i$, then either $j\in A$ in which case $A\in
S_{ij}$ or $j\nin A$ in which case $A\in S_{i/j}$. In either case,
$A\in \chi_{i}(S_{ij},S_{i/j}).\Box $  

In formally, this means that $S$ is a ``choice'' by player $i$ of
whether to join forces with player $j$ or to oppose him.

\begin{cor}\label{core}
If $S$ is a $n$-voter quota game ($n\geq 3$) then $S$ has weight
at most $n-2$. 
\end{cor}

{\em Proof:} By induction, and equations (\ref{part1}) and
(\ref{part2}).$\Box$ 

\cite[Lemma 2]{class} can be thought of as the converse of corollary
\ref{core}, since it uses methods similar to equations (\ref{part1})
and (\ref{part2}) to enumerate all quota games.
 
Note however that for some strong simple games, no pair of distinct
voters are comparable via influence. The only example with less than
seven vertices is the icosahedral game $I=S_{6,30}$: Associate each of six
voters with a pair of 
opposite vertices on an icosahedron. A coalition is winning if and
only if it contains a face. (See \S{A}.) 

Such games prevent us from generalizing the reasoning of corollary
\ref{core}. 

\appendix
\section{Classification of strong simple games}
In Table \ref{table}, we classify strong simple games
according to their height and weight. Table \ref{table} includes
all 30 isomorphism classes of strong simple games with up to six powerful
voters, as well as those with 
a transitive automorphism group and up to eight powerful voters, and
those of weight or height up to two.
Each isomorphism class is listed only once. 

Games $S$ are grouped according to their number $n$ of powerful voters
as indicated in the first column. The second
column provides the weights of a quota game if applicable, or
otherwise identifies the game. The remaining columns give
the weight $w$ and depth $d$ of the game along with some optimal
decompositions. Those case where our methods do not give (or are not
known to give) optimal
decompositions are marked $*$ (or $*?$). 
A question mark appears when optimality of $d$ or $w$ is uncertain. A
$\heartsuit $ 
appears in the first column if the game has a transitive automorphism
group. 
\begin{table}[htbp] \caption{Classification of Strong Simple Games}
\label{table}
{\scriptsize $$\scriptscriptstyle
\begin{array}{|l@{}r|l@{}r|l@{}r|}
\hline
n&S&w&\mbox{$m$-Median decomposition}&d&\mbox{$\chi$-Choice
decomposition}\\ 
\hline \hline
1\heartsuit&(1)_{1}&0&-&0&-\\ \hline
3\heartsuit&(111)_{2}
&	1&m((100)_{1},(010)_{1},(001)_{1})
&	1&\chi_{3}((100)_{1},(010)_{1})\\ \hline
4&(2111)_{3}
&	2&m((1110)_2,(1101)_2,(1011)_2) 
&	2&\chi_{4}((1000)_{1},(1110)_{2}) \\ \hline
5&(22111)_{4}
&	2& m((11100)_2, (11010)_2, (11001)_2) 
&	2& \chi_{5}((11100)_{2},(11010)_{2})  \\  
5&(31111)_{4}
&	2& m((10000)_1, (11100)_2, (10011)_2) 
&	2& *\chi_{1}((11100)_{2}, (10011)_{2}) \\
5&(32211)_{5}
&	2& m((11100)_2, (11010)_2, (10101)_2) 
&	3& \chi_{5}((21110)_{3}, (11100)_{2})\\
5\heartsuit&(11111)_{3}
&	3& m((20111)_3,(12011)_3,(01211)_3) 
&	3& \chi_{5}((01110)_{2}, (21110)_{3}) \\ \hline
6&(411111)_{5}
&	3& m((100000)_{1},(310111)_{4},(001000)_{1})
&	3& \chi_{6}((100000)_{1},(311110)_{4})\\
6&(522211)_{7}
&	3& m((100000)_{1},(320211)_{5},(001000)_{1})
&	3& \chi_{6}((311110)_{4},(211100)_{3})\\
6&(433111)_{7}
&	3& m((100000)_{1},(130111)_{4},(001000)_{1})
&	3& \chi_{3}((130111)_{4},(100000)_{1}) \\
6&(533211)_{8}
&	3& m((100000)_{1},(230211)_{5},(001000)_{1})
&	4?&\chi_{6}((211100)_{3},(322110)_{5})\\
6&(422111)_{6}
&	2&m((100000)_1, (111000)_2, (000111)_2)
&	2&*\chi_{1}((111000)_{2},(000111)_{2})\\ 
6&(543221)_{9}
&	3&m((100000)_{1},(120110)_{3},(013111)_{4})
&	3&\chi_{3}((120110)_{3},(310111)_{4})\\
6&(332111)_{6}
&	2&m((111000)_2, (110100)_2, (001011)_2) 
&	3&\chi_{6}((111000)_{2},(221110)_{4}) \\
6&(432211)_{7}
&	2&m((111000)_2, (110100)_2, (100011)_2) 
&	3&\chi_{6}((211100)_{3},(221110)_{4})\\ 
6&(542222)_{9}
&	3&m((100000)_{1},(220111)_{4},(013111)_{4})
&	3&*?\chi_{1}((220111)_{4},(013111)_{4})\\
6&(332221)_{7}
&	3& m((301111)_{4},(120110)_{3},(012110)_{3})
&	3& \chi_{2}((301111)_{4},(001110)_{2})\\
6&(222111)_{5}
&	3& m((301111)_{4},(130111)_{4},(013111)_{4})
&	3& \chi_{2}((301111)_{4},(002111)_{3})\\
6&(322211)_{6}
&	3& m((301111)_{4}, (230211)_{5}, (013111)_{4})
&	3& \chi_{2}((301111)_{4},(102211)_{4})\\
6&(211111)_{4}
&	3& m((301111)_{4}, (220111)_{4}, (013111)_{4})
&	3&\chi_{6}((211100)_{3},(221110)_{4})\\
6&S_{6,21}
&	3&m((100000)_{1},(120101)_{3},(012110)_{3})
&	3&*?\chi_{1}((120101)_{3},(012110)_{3})\\
6&S_{6,22}
&	3&m((100000)_{1},(221110)_{4},(010112)_{3})
&	3&*?\chi_{1}((221110)_{4},(010112)_{3})\\
6&S_{6,23}
&	3&m((031111)_{4},(102110)_{3},(310111)_{4})
&	3&\chi_{2}((200111)_{3},(102110)_{3})\\
6&S_{6,24}
&	3&m((021110)_{3},(102101)_{3},(210110)_{3})
&	3&\chi_{6}((122110)_{4},(211210)_{4})\\
6&S_{6,25}
&	3&m((031111)_{4},(102101)_{3},(310111)_{4})
&	3&\chi_{6}((122110)_{4},(221110)_{4})\\
6&S_{6,26}
&	3&m((021110)_{3},(203112)_{5},(210110)_{3})
&	4?&\chi_{6}((122110)_{4},(311220)_{5})\\
6&S_{6,27}
&	2&m((000001)_1, (111000)_2, (001110)_2)
&	2&\chi_{6}((111000)_{2},(001110)_{2})\\
6&S_{6,28}
&	2& m((110001)_2, (101010)_2, (011100)_2)
&	3& \chi_{3}((110001)_{2},(110221)_{4})\\
6&S_{6,29}
&	3&m((200111)_{3},(020111)_{3},(002111)_{3})
&	3&\chi_{1}((020111)_{3},(002111)_{3})\\
6\heartsuit&I=S_{6,30}
&	3&m((302112)_{5}, (230121)_{5}, (023211)_{5})
&	5?& \mbox{*\footnotemark{ }}m(S_{6,28},(111220)_{4})
\\ \hline
7&B_{2}
&	2&m((1110000)_{2}, (0001110)_{2}, (0000001)_{1}) 
&	2&*\chi_{7}((1110000)_{2},(0001110)_{2})\\
7&
&	2&m((1110000)_{2}, (0001110)_{2}, (0000111)_{2})
&	3&\chi_{1}((0011221)_{4},(0101221)_{4})\\
7&
&	2&m((1110000)_{2}, (1001100)_{2}, (1000011)_{2})
&	3&\chi_{2}(S_{6,27},(3001111)_{4}) \\
7&
&	2&m((1110000)_{2}, (0011100)_{2}, (0000111)_{2})
&	3&\chi_{7}(S_{6,27},S'_{6,27}) \\ 
7\heartsuit&\mbox{Fano}
&	4?&m(S_{6,22},S'_{6,22},S''_{6,22})
&	5?&\mbox{*\footnotemark{ }}
	\chi_{7}(S_{6,22},\chi_{7}(S'''_{6,22},S''''_{6,22})) \\
7\heartsuit&(1111111)_{4}
&	3&m((2011111)_{4},(1201111)_{4},(0121111)_{4})
&	4?&\chi_{7}((1111100)_{3},(2111110)_{4})\\ \hline
8&&2&m\left(\begin{array}{c}
(11100000)_{2},(00011100)_{2},\\
(00000111)_{2}
\end{array} \right)&3&\chi_{8}(S_{6,27},S'_{6,27})\\ \hline
9\heartsuit&Dem_{3}^{2}&2&m\left(\begin{array}{c}
(111000000)_{2},(000111000)_{2},\\
(000000111)_{2}
\end{array} \right)
&3&\chi_{9}(B_{2},B'_{2}) \\
\hline
\end{array}$$}
\end{table}
\addtocounter{footnote}{-1}
\footnotetext{Proposition \ref{decompose} does not apply to games
whose influence relation is trivial such as the Icosahedral game $I$.}
\addtocounter{footnote}{1}
\footnotetext{The players in the Fano game are
identified with points in the projective plane over the field ${\bf
Z}_{2}$. A coalition wins if it contains a line. Proposition
\ref{decompose} does not apply.}

\end{document}